\documentclass[11pt]{article}
\usepackage{amsmath,amssymb}

\newtheorem{propo}{{\bf Proposition}}[section]
\newtheorem{coro}[propo]{{\bf Corollary}}
\newtheorem{lemma}[propo]{{\bf Lemma}} \newtheorem{theor}[propo]{{\bf
Theorem}}

\begin{document}

\vspace*{1.0in}

\begin{center}FURTHER RESULTS ON ELEMENTARY LIE ALGEBRAS AND LIE A-ALGEBRAS

\end{center}
\bigskip

\centerline {David A. Towers} \centerline {Department of
Mathematics, Lancaster University} \centerline {Lancaster LA1 4YF,
England}

\centerline {and}

\centerline {Vicente R. Varea \footnote[1]{Supported by DGI Grant
BFM2000-1049-C02-01}} \centerline {Department of Mathematics,
University of Zaragoza} \centerline {Zaragoza, 50009 Spain}

\bigskip

{\bf Abstract} A finite-dimensional Lie algebra $L$ over a field $F$ of characteristic zero is called elementary if each of its subalgebras has trivial Frattini ideal; it is an $A$-algebra if every nilpotent subalgebra is abelian. This paper is a continuation of the study of these algebras initiated by the authors in \cite{tow-var}. If we denote by $\mathcal{A}$, $\mathcal{G}$, $\mathcal{E}$, $\mathcal{L}$, $\Phi$ the classes of $A$-algebras, almost algebraic algebras, $E$-algebras, elementary algebras and $\phi$-free algebras respectively, then it is shown that:
\[ \mathcal{L} \subset \Phi \subset \mathcal{G}, \hspace{.2cm} \mathcal{L} \subset \mathcal{A} \subset \mathcal{E} \hspace{.2cm} \hbox{and} \hspace{.2cm} \mathcal{G} \cap \mathcal{A} = \mathcal{L}.
\]
It is also shown that if $L$ is a semisimple Lie algebra all of whose minimal parabolic subalgebras are $\phi$-free then $L$ is an $A$-algebra, and hence elementary. This requires a number of quite delicate properties of parabolic subalgebras. Finally characterisations are given of $E$-algebras and of Lie algebras all of whose proper subalgebras are elementary.
\medskip


\bigskip 

\section{Introduction}
This paper is a continuation of the study initiated by the authors in \cite{tow-var}. Throughout $L$ will denote a finite-dimensional Lie algebra over a field $F$. The Frattini ideal of $L$, $\phi(L)$, is the largest ideal of $L$ contained in all maximal subalgebras of $L$. The Lie algebra $L$ is called {\em $\phi$-free} if $\phi(L) = 0$, and {\em elementary} if $\phi(B)=0$ for every subalgebra $B$ of $L$. Lie algebras all of whose nilpotent subalgebras are abelian are called {\em $A$-algebras}; Lie algebras $L$ such that $\phi(B)\leq \phi(L)$ for all subalgebras $B$ of $L$ are called {\em $E$-algebras}. We are seeking to determine properties of, and inter-relationships between, these three classes of algebras. 
A linear Lie algebra
$L\leq {\rm gl}(V)$ is {\it almost algebraic} if $L$ contains the
nilpotent and semisimple Jordan components of its elements; an abstract Lie algebra
$L$ is then called almost algebraic if ${\rm ad}L\leq {\rm gl}(L)$ is
almost algebraic. 
\par
Throughout sections two and three $F$ is assumed to have characteristic zero. In section 2 we show first that $L$ is $\phi$-free if and only if it is almost algebraic and its nilradical is abelian. It follows from this that if $L$ is almost algebraic then $\phi(L) = N^2$, where $N$ is the nilradical of $L$. Zhao and Lu proved in \cite{zl} that every almost-algebraic $A$-algebra is elementary, whenever the ground field is algebraically closed of characteristic zero. We generalise this by showing that $L$ is an almost-algebraic $A$-algebra if and only if it is elementary (and without the assumption of an algebraically closed field). From this we deduce that if $L$ is a Lie $A$-algebra with a $\phi$-free solvable radical, then $L$ is elementary. The final result in this section is that if $L$ is a Lie $A$-algebra then it is an $E$-algebra.
\par
A subalgebra $P$ of $L$ is called {\it parabolic} if $P \otimes_F \overline{F}$ contains a Borel subalgebra (that is, a maximal solvable subalgebra) of $L \otimes_F \overline{F}$, where $\overline{F}$ is the algebraic closure of $F$. The main purpose of section 3 is to prove two results. The first is that if $P$ is a minimal parabolic subalgebra of $L$ then the following are equivalent: $P$ is $\phi$-free; the nilradical of $P$ is abelian; and $P$ is elementary. The second is that if $L$ is a semisimple Lie algebra all of whose minimal parabolic subalgebras are $\phi$-free then $L$ is an $A$-algebra, and hence elementary. In order to establish these results we need a number of properties of parabolic subalgebras. Some of these may be known, but we know of no references to them other than that V\"{o}lklein in \cite{volk} shows that if $L$ is semisimple then the minimal parabolic subalgebras of $L$ are the idealisers of the maximal nil subalgebras of $L$. He uses, however, the canonical correspondence between the connected subgroups of the identity component of the automorphism group of $L$ and the algebraic subalgebras of $L$. Our proofs are based entirely on internal properties of $L$ itself and so, we believe, are of interest in themselves. We also include these proofs for the convenience of the reader.
\par
In the final section we look more closely at Lie $E$-algebras. In particular we give a characterisation of them over a field of characteristic zero. We also characterise Lie algebras all of whose proper subalgebras are elementary. These generalise Theorems 4.3 and 4.7 and Corollaries 4.4 and 4.5 of \cite{tow} by removing the requirement that the underlying field be algebraically closed. 
\par
We will denote vector space direct sums by $\oplus$ and semidirect products 
by $\rtimes$. If $A$ is a subalgebra of $B$ we will write $A \leq B$, whereas $A < B$ will mean that $A$ is a proper subalgebra of $B$. The (solvable) {\em radical} (resp. {\em nilradical}) of $L$ will be denoted by $R(L)$ (resp. $N(L)$), whilst $Asoc(L)$ will denote the sum (necessarily direct) of the minimal abelian ideals of $L$.
\bigskip

\section{Preliminary results}
First we reveal the relationship between almost-algebraic and $\phi$-free Lie algebras.
\bigskip

\begin{propo}\label{p:phifree} Let $L$ be a Lie algebra over a field of characteristic zero.
\begin{itemize} 
\item[(i)] If $L$ is $\phi$-free, then $L$ is
almost algebraic. 
\item[(ii)] Let $L$ be almost algebraic. Then $L$ is
$\phi$-free if and only if its nilradical is abelian.
\end{itemize} 
\end{propo} 
\medskip
{\it Proof.} (i) Let $L$ be
$\phi$-free. By \cite{frat} we have that $L = N(L) + S$ where $S$
is a subalgebra of $L$ such that $\hbox{ad}_L S$ is completely
reducible and $N(L) ={\rm Asoc} (L)$. From \cite[Theorem 2.2]{ab}
it follows that $L$ is almost algebraic.

(ii) By \cite[Theorem 2.2]{ab} we have
that $L = N(L) + S$ where $S$ is a subalgebra of $L$ such that
$\hbox{ad}_L S$ is completely reducible. Now, assume that $N(L)$ is abelian.
Then we have that $N(L) =$ Asoc $L$. So $L$ is $\phi$-free by \cite{frat}.
The converse follows from \cite{frat}.
\bigskip

\begin{coro}\label{c:aalg} Let $L$ be an almost-algebraic Lie algebra over a field of
characteristic zero. Then $\phi(L) = N^2$, where $N$ is the nilradical of $L$.
\end{coro}
\medskip
{\it Proof.} Clearly $N^2$ is almost algebraic and so $L/N^2$ is
almost algebraic, by \cite[Lemma 4.1]{ab}. Moreover, $N(L/N^2)$ is
abelian and hence $\phi$-free by Proposition \ref{p:phifree}. It
follows that $\phi(L) \subseteq N^2$. The reverse inclusion is
given by \cite[Theorem 6.5]{frat}.
\bigskip

If $B$ is a subalgebra of $L$ we define the {\em idealiser} of $B$ in $L$ to be $I_L(B) = \{x \in L : [x,B] \leq B \}$. Next we show that a Lie algebra is an almost-algebraic $A$-algebra if and only if it is elementary, thereby generalising the result of Zhao and Lu proved in \cite{zl}. First we need the following result.
\bigskip

\begin{propo}\label{p:aaalg} Let $L$ be an almost-algebraic Lie algebra over a field of characteristic zero. If every almost-algebraic subalgebra of $L$
is $\phi$-free, then $L$ is elementary. 
\end{propo}
\medskip
{\it Proof.} Let $B$ be a subalgebra of $L$. Then the idealiser of
$B$ in $L$, $I_L(B)$, is almost algebraic, by \cite[Theorem
2.3]{ab}. By our hypothesis, we have $\phi(I_L(B)) = 0$. But now
$B$ is an ideal of $I_L(B)$, so $\phi(B) \leq \phi(I_L(B)) $, by
\cite[Corollary 4.2]{frat}. Hence $L$ is elementary.
\bigskip

\begin{theor}\label{p:aalg}
Let $L$ be a Lie algebra over a field of characteristic zero.
Then $L$ is an almost-algebraic $A$-algebra if and only if it is elementary.
\end{theor}
\medskip
{\it Proof.} ($\Rightarrow$) Let $L$ be an almost-algebraic
$A$-algebra and let $B$ be an almost-algebraic subalgebra of $L$.
Then by Corollary \ref{c:aalg} it follows that $\phi(B) = N^2$,
where $N$ is the nilradical of $B$. Therefore, $\phi(B)=0$ since
$N$ is abelian. Hence $L$ is elementary by Proposition
\ref{p:aaalg}.
\par
($\Leftarrow$) This follows from \cite[Corollary 4.7]{tow-var}.
\bigskip

\begin{coro}\label{l:rad} Let $L$ be a  Lie $A$-algebra over a
field of characteristic zero. If $R(L)$ is $\phi$-free, then $L$
is elementary. 
\end{coro} 
\medskip 
{\it Proof.} Assume that $R(L)$ is $\phi$-free. By Proposition \ref{p:phifree} we have that $R(L)$ is almost algebraic. Then, from \cite[Corollary 3.1]{ab} it
follows that $L$ is also almost algebraic. So, $L$ is elementary
by Theorem \ref{p:aalg}. \bigskip

\begin{coro}\label{c:split}
Let $L$ be an almost-algebraic Lie $A$-algebra over a field of characteristic zero.
Then $L$ splits over each of its ideals.
\end{coro}
\medskip
{\it Proof.} This follows from Theorem \ref{p:aalg} and
\cite[Lemma 2.3]{tow}.
\bigskip

Finally we have that Lie $A$-algebras are necessarily $E$-algebras.
\bigskip

\begin{propo}\label{p:aealg} Let $L$ be a Lie $A$-algebra over a
field of characteristic zero. Then $L$ is an $E$-algebra.
\end{propo} 
\medskip 
{\it Proof.} We have that $L/\phi(L)$ is an $A$-algebra, by \cite[Lemma 1]{prem}. On the other hand, we have that $L/\phi(L)$ is $\phi$-free and so almost algebraic by Proposition \ref{p:phifree}. Then $L/\phi(L)$ is elementary, by Theorem \ref{p:aalg}. So, $L$ is an $E$-algebra, by \cite[Proposition 2]{stit}. 
\bigskip

Denote by $\mathcal{A}$, $\mathcal{G}$, $\mathcal{E}$, $\mathcal{L}$, $\Phi$ the classes of $A$-algebras, almost algebraic algebras, $E$-algebras, elementary algebras and $\phi$-free algebras respectively. Then, to summarise, what we have shown is the following:
\[ \mathcal{L} \subset \Phi \subset \mathcal{G}, \hspace{.2cm} \mathcal{L} \subset \mathcal{A} \subset \mathcal{E} \hspace{.2cm} \hbox{and} \hspace{.2cm} \mathcal{G} \cap \mathcal{A} = \mathcal{L}.
\]

\section{Parabolic subalgebras}

Throughout this section $L$ denotes a (non-zero) semisimple Lie algebra over a field $F$ of characteristic zero. We denote by $\overline{F}$ the algebraic closure of $F$, and write $\overline{S} = S \otimes_F \overline{F}$
for each subspace $S$ of $L$. For each subalgebra $S$ of $L$, let
$U(S)$ denote the set of ad-nilpotent elements in the solvable
radical, $R(S)$, of $S$. Our main objective in this section is to show that in order to check whether a Lie algebra is elementary it suffices to look at its minimal parabolic subalgebras. First we have some properties of the nilradical of a parabolic subalgebra.
\bigskip

\begin{propo}\label{p:nilradical} Let $P$ be a
parabolic subalgebra of $L$. Then 
\begin{itemize} 
\item[(i)] $N(P)$ is nil in $L$; 
\item[(ii)] if $P< Q< L$, then $N(Q)< N(P)$; and 
\item[(iii)] $P=I_L(N(P))$.
\end{itemize}
\end{propo}\medskip

{\it Proof.} Since $\overline{N(P)}$ is the nilradical of
$\overline{P}$ and $\overline{I_L(N(P))} =
I_{\overline{L}}(\overline{N(P)})$ (see \cite[pages 42 and 36]{baki}) we may assume that $F$ is
algebraically closed. Now $P$ is conjugate to a standard parabolic
subalgebra and so we can assume it is of the following form. Let $H$
be a Cartan subalgebra of $L$ and let $\Delta$ be the set of roots
corresponding to $H$. Then \[ P = H \oplus
\Sigma_{\alpha\in\Delta^+}L_{\alpha} \oplus \Sigma_{\alpha \in
\Omega_1} L_{\alpha}, \] where $\Omega_1\subseteq\Delta^-$. Let
$\Omega_1^{\prime}=\{\alpha\in\Delta^+\mid
-\alpha\not\in\Omega_1\}$. Then we have that \[ N(P) =
\Sigma_{\alpha \in \Omega_1^{\prime}} L_{\alpha}\] It follows that
$N(P)$ is nil in $L$ and (i) is proved.
\par
To prove (ii), suppose that $P < Q < L$. Then $Q$ is also parabolic
and so has the same form as $P$ but with $\Omega_1$ replaced by
$\Omega_2$ where $\Omega_1 \subset \Omega_2$. We have that
$\Omega_2^{\prime}\subset \Omega_1^{\prime}$ and hence $N(Q) <
N(P)$.
\par
To prove (iii), put $Q=I_L(N(P))$ and suppose that $P < Q$. By (ii)
we have $N(Q)< N(P)$. On the other hand, as $N(P)$ is a nilpotent
ideal of $Q$ we have $N(P)\leq N(Q)$, which is a contradiction.
Now the proof is complete.
\bigskip

The {\em centre} of $L$ is the set $Z(L) = \{ x \in L : [x,L] = 0 \}$. A subalgebra $T$ of $L$ is said to be a {\em toral} subalgebra of
$L$ if $T$ is abelian and ${\rm ad}_Lt$ is semisimple for every
$t\in T$. Next we need that if $P$ is an algebraic subalgebra of $L$ then $U(P)$ behaves well under field extension.
\bigskip

\begin{lemma}\label{l:nil} Let $U$ be a nil subalgebra of $L$. Then $\overline
U$ is also a nil subalgebra of $\overline L$
\end{lemma}

{\it Proof.} We have that $\overline U$ is a nilpotent algebraic
subalgebra of $\overline L$, since these properties are inherited from $U$ (see
\cite[p.181]{Chevalley}). So, $\overline U=U(\overline U)\oplus
T$, where $T$ is a toral subalgebra of $\overline L$ and
$[T,U(\overline U)]=0$, by \cite[Theorem 4]{chev}. We have $T\leq Z(\overline
U)=\overline{Z(U)}$ (see \cite[page 36]{baki}). As $Z(U)$ consists of commuting ad-nilpotent elements, we have that $\overline{Z(U)}$ is a nil subalgebra of
$\overline L$. This yields that $T$ is both nil and toral in $\overline
L$ and so $T\leq Z(\overline L)=0$. This yields that $\overline
U=U(\overline U)$ and so $\overline U$ is nil in $\overline
L$.
\bigskip

\begin{propo}\label{p:u} Let $P$ be an algebraic subalgebra of $L$. Then
$\overline{U(P)}=U(\overline P)$
\end{propo}

{\it Proof.} Since $P$ is algebraic, we have that $P=U(P)\oplus
M$, where $M=S\oplus Z(M)$ and $S$ is semisimple ($M$ is a Levi
factor of $P$), by \cite[Theorem 4]{chev}. It follows that
$\overline P=\overline {U(P)}\oplus\overline  M$ and $\overline
M={\overline S}\oplus Z(\overline M)$. On the other hand, from
Lemma \ref{l:nil}, it follows that $\overline{U(P)}$ is a nil
ideal of $\overline L$. So, $\overline {U(P)}\leq U(\overline P)$.
This yields that $\overline {U(P)}= U(\overline P)$.
\bigskip

The following results are concerned with relationships between parabolic subalgebras and certain nil subalgebras.
\bigskip

\begin{lemma}\label{l:p} Let $U$ be a nil subalgebra of $L$ and put $P=I_L(U)$.
If $U=U(P)$, then $P$ is parabolic.
\end{lemma}
\medskip
{\it Proof.} This follows from Proposition \ref{p:u} and
\cite[Theorem 29.8.1]{ty}

\begin{propo}\label{p:p} Let $U$ be a nil subalgebra of $L$ and put $P=I_L(U)$. Then there is a parabolic subalgebra $Q$ of $L$ satisfying:
\begin{itemize}
\item[(i)] $U\leq U(Q)$;
\item[(ii)] $P \leq Q$; and
\item[(iii)] $U(P)\leq U(Q)$.
\end{itemize}
\end{propo}
\medskip
{\it Proof.} Put $U_0=U$ and $Q_1=I_L(U_0)$. Define inductively
the two sequences $\{U_i\}_{i\geq 0}$, $\{Q_i\}_{i\geq 1}$ by
$Q_i=I_L(U_{i-1})$ $U_i=U(Q_i)$. Then
$Q_i=I_L(U_{i-1})=I_L(U(Q_{i-1}))\geq Q_{i-1}$, so these sequences
are increasing. This yields that there is an integer $j$ such that
$Q_j=Q_{j+1}$; that is $Q_j=I_L(U_j)$. It follows from Lemma \ref
{l:p} that $Q=Q_j$ is a parabolic subalgebra of $L$. We have that
$U\leq U_j=U(Q)$, giving (i), and $P = I_L(U)=Q_1\leq Q$, giving (ii).
\par
Finally, $Q_i=I_L(U_{i-1})$ so $U_{i-1}$ is a nil ideal of $Q_i$
and $U_{i-1}\leq U(Q _i)$, whence $U_1=U(P)\leq U(Q)$, giving (iii).

\begin{propo}\label{p:minimal} 
\begin{itemize}
\item[(i)] If $U$ is a maximal nil subalgebra of $L$, then $I_L(U)$ is a minimal
parabolic subalgebra of $L$.
\item[(ii)] If $P$ is a minimal parabolic subalgebra of $L$, then $U(P)$ is a maximal nil subalgebra of $L$
\end{itemize}
\end{propo}
\medskip
{\it Proof.} (i): Let $U$ be a maximal nil subalgebra of $L$. We
have that $U$ is a nil ideal of $I_L(U)$ and so $U\leq U(I_L(U)$.
By the maximality of $U$, we must have that $U=U(I_L(U))$. From
Lemma \ref{l:p} it follows that $I_L(U)$ is parabolic. Now let $Q$
be a parabolic subalgebra with $Q \leq I_L(U)$. By Proposition
\ref{p:nilradical} we have that $U = N(I_L(U)) \leq
N(Q) = U(Q)$. By the maximality of $U$ it follows that $U=N(Q)$.
This yields that $I_L(U)=Q$ by Proposition \ref{p:nilradical} again, and
so $I_L(U)$ is a minimal parabolic subalgebra of $L$.
\par 
(ii): Let $P$ be a minimal parabolic subalgebra of $L$. Suppose that
there is a nil subalgebra $V$ of $L$ such that $U(P)< V$. By
Engel's Theorem we have that $U(P)<V\cap I_L(U(P))$. As $P$ is
parabolic, $I_L(U(P))=P$. Put $W=V\cap P$. Let $M$ be a Levi
factor of $P$; so that $P=U(P)\oplus M$, $M=Z(M)\oplus S$ where
$Z(M)$ is toral in $L$ and $S$ is a semisimple subalgebra of $L$.
\par
We see that $W=U(P)\oplus (M\cap V)$. Since $Z(M)$ is toral it
follows that $M\cap V=S\cap V$. So, $S\cap V$ is a non-trivial nil
subalgebra of the semisimple Lie algebra $S$. From Proposition
\ref{p:p} there is a parabolic subalgebra $Q$ of $S$, $Q\not=S$,
containing $S\cap V$. Let $B$ a Borel subalgebra of $\overline{S}$
contained in $\overline Q$. Then, we have that $R(\overline P)+B$
is a maximal solvable subalgebra of $\overline P$. Since
$\overline P$ is parabolic in $\overline L$, it follows that
$R(\overline P)+B$ is a Borel subalgebra of $\overline L$. We have
$\overline{R(P)+Q}=R(\overline P)+\overline Q\geq R(\overline
P)+B$. Therefore, $R(P)+Q$ is a parabolic subalgebra of $L$
contained in $P$, which contradicts the minimality of $P$. Hence
$U(P)$ is a maximal nil subalgebra of
$L$.
\bigskip

A Lie algebra $L$ is said to be {\em ad-semisimple} if ${\rm
ad}x$ is semisimple for every $x\in L$. Then we have the following criterion for a parabolic subalgebra to be minimal.
\bigskip

\begin{coro}\label{c:min} A parabolic subalgebra $P$ of $L$ is
minimal if and only if $P/U(P)$ is ad-semisimple.
\end{coro}
\medskip
{\it Proof.} Let $M$ be a Levi factor of $P$. Let us first suppose
that $P$ is minimal. Then, by Proposition \ref{p:minimal} it
follows that $U(P)$ is a maximal nil subalgebra of $L$ and so $M$
is ad-semisimple. Now assume that $P/U(P)$ is ad-semisimple, so
that $M$ is ad-semisimple. As $Z(L)=0$, we see that $U(P)$ is a
maximal nil subalgebra of $P$. Since $P=I_L(P)$, it follows from Engel's
Theorem that $U(P)$ is maximal nil subalgebra of $L$.
\bigskip

We now have the results that we need to show the role played by the minimal parabolic subalgebras in determining whether or not $L$ is elementary.
\bigskip

\begin{lemma}\label{l:var}
Let $P$ be a minimal parabolic subalgebra of $L$. Then the
following are equivalent
\begin{itemize}
\item[(i)] $P$ is $\phi$-free;
\item[(ii)] $N(P)$ is abelian; and
\item[(iii)] $P$ is elementary.
\end{itemize}
\end{lemma}
\medskip
{\it Proof.} (i)$\Rightarrow$(ii): Since $P$ is algebraic, it follows from
Corollary \ref{c:aalg} that $N(P)^2=\phi(P)=0$.
\par
(ii)$\Rightarrow$(iii): Since $N(R(P))=N(P)$ and $R(P)$ is algebraic,
it follows from Proposition \ref{p:phifree} that $R(P)$ is
$\phi$-free. By \cite[Theorem 2.5]{tow-var} we have that $R(P)$ is
elementary. On the other hand, since $U(P)\leq R(P)$, it follows from
Corollary \ref{c:min} that $P/R(P)$ is ad-semisimple.
But \cite[Proposition 4.4]{tow-var} now implies that $P$ is
elementary.
\par
(iii)$\Rightarrow$(i): This is trivial.
\bigskip

\begin{theor}\label{p:elem}
Let $L$ be a semisimple Lie algebra over a field of characteristic
zero, and suppose that all minimal parabolic subalgebras of $L$
are $\phi$-free. Then $L$ is an $A$-algebra, and hence elementary.
\end{theor}
\medskip
{\it Proof.} First we show that $L$ is an $A$-algebra. Let $U$ be
a nilpotent subalgebra of $L$. Then $U$ is contained in a maximal
solvable subalgebra $\Gamma$ of $L$. As $\Gamma$ is algebraic we
can write $\Gamma = U(\Gamma) \dot{+} T$ where $T$ is a toral
subalgebra. Let $N$ be a maximal nil subalgebra of $L$ containing
$U(\Gamma)$, and let $P$ be the idealiser in $L$ of $N$. By
Proposition \ref{p:minimal}, $P$ is minimal parabolic and hence
$\phi$-free. It follows from Lemma \ref{l:var} that $P$ is
elementary, and hence an $A$-algebra by Proposition \ref{p:aalg}.
This yields that $U(\Gamma)$ is abelian. But then $\Gamma$ is
$\phi$-free and hence elementary, by \cite[Theorem 2.5]{tow-var}.
It follows that $U$ is abelian and that $L$ is an $A$-algebra.
\par
By \ref{p:aalg} it follows that $L$ is elementary.
\bigskip

\section{$E$-algebras}
First we have an easy strengthening of Corollary 2.2 of \cite{tow-var}.
\bigskip

\begin{propo}\label{p:esolv}
Let $L$ be a solvable Lie algebra over a perfect field. Then $L$
is an $E$-algebra if and only if $L$ is strongly solvable.
\end{propo}
\medskip
{\it Proof.} This follows from \cite[Corollary 2.2]{tow-var} and
\cite{stit}.
\bigskip

Next we have the following versions of Theorem 4.3 and Corollaries 4.4 and 4.5 of \cite{tow} with the assumption that the underlying field be algebraically closed removed.
\bigskip

\begin{theor}\label{t:egen}
Let $L$ be a Lie algebra over a field of characteristic zero. Then
$L$ is an $E$-algebra if and only if one of the following holds:
\begin{itemize}
\item[(i)] $L$ is solvable; 
\item[(ii)] $L$ is elementary and semisimple; or
\item[(iii)] $L = R \oplus S$,
where $R = R(L)$ is the radical of $L$, $S = S_1 \oplus S_2$,
$S_1$ is an ad-semisimple ideal of $S$, $S_2$ is an elementary and semisimple ideal of $S$, and $S_2 R \leq \phi(L)$.
\end{itemize}
\end{theor}
\medskip
{\it Proof.} ($\Rightarrow$): Let $L$ be an $E$-algebra and
suppose that $L$ is not solvable or semisimple. By Levi's Theorem
we can write $L = R \oplus S$
 where $S$ is a semisimple subalgebra of $L$.
 If $U$ is a subalgebra of $L$ then write $\overline{U}$ for its image in
$L/\phi(L)$ under the natural homomorphism.
 We have $\overline{L} = \overline{R} \oplus \overline{S}$ and $\overline{S} \cong S$ is
elementary and semisimple.
 Put $S = S_1 \oplus S_2$ where $S_2$ is the largest semisimple ideal of
$L$.
 If $\overline{S_1} \neq 0$, then $\overline{S_1} \cong S_1$ is ad-semisimple, as
in \cite[Theorem 4.6]{tow-var}.
 Clearly $S_2 R \leq \phi(L)$, so (iii) holds.

($\Leftarrow$): If (i) holds then $L$ is strongly solvable (since
the ground field has characteristic 0) and the result follows from
Proposition \ref{p:esolv} above. If (ii) holds the result is
clear. So suppose that (iii) holds. Then $\overline{L} = \overline{L_1}
\oplus \overline{S_2}$ where $\overline{L_1} = \overline{R} \oplus \overline{S_1}$.
Now $\overline{L_1}/ \overline{R} \cong \overline{S_1} \cong S_1$ is
ad-semisimple and $\overline{R}$ is elementary, so $\overline{L_1}$ is
elementary, by \cite[Proposition 4.4]{tow-var}. It follows that
$\overline{L}$ is elementary, and hence that $L$ is an $E$-algebra.
\bigskip

\begin{coro}\label{c:enilp}
Let $L$ be a Lie algebra over a field of characteristic zero, and
suppose that the radical of $L$ is nilpotent. Then $L$is an
$E$-algebra if and only if one of the following holds:
\begin{itemize}
\item[(i)] $L$ is nilpotent;
\item[(ii)] $L$ is elementary and semisimple; or
\item[(iii)] $L = (R \oplus S_1) \oplus S_2$, where $R = R(L)$ is the
radical of $L$, $S_1$ is an ad-semisimple subalgebra of $L$, and $S_2$ is an elementary and semisimple ideal of $L$.
\end{itemize}
\end{coro}
\medskip
{\it Proof.} This follows as in \cite[Corollary 4.4]{tow}.
\bigskip

\begin{coro}\label{c:eperf}
Let $L$ be a perfect Lie algebra (i.e., $L = L^2$) over a field of
characteristic zero. Then $L$ is an $E$-algebra if and only if $L$
is elementary and semisimple.
\end{coro}
\bigskip

Finally we consider non-elementary Lie algebras all of whose proper subalgebras are elementary. We call such algebras {\em minimal non-elementary} Lie algebras. The following extends Theorem 4.7 of \cite{tow}.
\bigskip

\begin{theor}\label{t:alel} Let $L$ be a Lie algebra over a field
of characteristic zero. Then $L$ is a minimal non-elementary Lie
algebra if and only if 
\begin{itemize} 
\item[(i)] $L = L^2 \rtimes
Fx$, where $L^2$ is abelian and $0 \neq \phi(L) = Asoc L$ is the
biggest ideal of $L$ properly contained in $L^2$, or 
\item[(ii)] $L$ is the three-dimensional Heisenberg algebra. 
\end{itemize}
\end{theor} 
\medskip 
{\it Proof.} ($\Rightarrow$) First note that
$L$ must be an $E$-algebra. Suppose that $L$ is not solvable. Then
$R = R(L)$ is elementary and so almost algebraic by
\cite[Proposition 4.1]{tow-var}. Hence $L$ is almost algebraic, by
\cite[Corollary 3.1]{ab}. Moreover, the nilradical is elementary
and so abelian. It follows that $L$ is $\phi$-free, by Proposition
\ref{p:phifree}, and so elementary - a contradiction. This yields
that $L$ is solvable. 
\par 
Suppose that $L$ is not nilpotent. Then
$\phi(L) \neq L^2$ (see, for example, \cite[section 5]{frat}), so
there is a maximal subalgebra $M$ of $L$ such that $L = L^2 + M$.
Choosing $B$ to be a subalgebra minimal with respect to the property that $L =
L^2 + B$ we have $L^2 \cap B \leq \phi(B) = 0$ (see \cite[Lemma
7.1]{frat}), so $L = L^2 \oplus B$ and $B$ is abelian. Moreover, $L^2$ is nilpotent
and elementary, and so abelian. 
\par 
Suppose that dim$B > 1$. Let $K$ be a maximal subalgebra of
$B$. Then $M = L^2 + K$ is a maximal subalgebra of $L$, so $M$
is elementary and hence $\phi$-free. It follows from \cite[Theorem
7.4]{frat} that $L^2$ is completely reducible as a $K$- module,
and hence that each element of $K$ acts semisimply on $L^2$
(\cite[Theorem 10, page 81]{jac}). But every element of $B$ is
contained in a maximal subalgebra of $B$ and so acts semisimply on
$L^2$. This yields that $L^2$ is a completely reducible
$B$-module, whence $L^2 \leq$ Asoc $L$ and $L$ splits over Asoc
$L$. But then $\phi(L) = 0$, by \cite[Theorem 7.3]{frat}, and $L$
is elementary, a contradiction. Thus dim$B = 1$, so put $B = Fx$.
\par 
So we now have $L = L^2 \rtimes Fx$ where $L^2$ is abelian.
Let $C$ be an ideal of $L$ with Asoc $L < C < L^2$ and put $D = C
+ Fx$. Then $D \neq L$ so $\phi(D) = 0$, giving Asoc $D = N(D)
\geq C$. But every minimal ideal of $D$ is inside $L^2$ and
invariant under ad$x$ and so is an ideal of $L$. It follows that
$C =$ Asoc $D \leq$ Asoc $L$, a contradiction. Hence Asoc $L$ is
the biggest ideal of $L$ properly contained in $L^2$. We must have
Asoc $L \leq \phi(L)$, since otherwise $L$ splits over Asoc $L$,
as in paragraph two above. As $\phi(L) \neq L^2$ this means that
$\phi(L) =$ Asoc $L$ and we have case (i). 
\par 
Suppose now that $L$ is nilpotent. Then $L$ is not abelian, so dim $L \geq 3$ 
and $L$ has a chain of ideals 
\[ 0 = L_0 < L_1 < ... < L_n = L, 
\]
where dim $L_i = i$ and $LL_i \leq L_{i-1}$ for $1 \leq i \leq n$.
Let $L_1 = Fz$, $L_2 = Fy + Fz$, and let $x$ be any element of
$L$. If $Fx + Fy + Fz$ is abelian for every $x \in L$ then $L$ is
abelian, a contradiction. Hence $L$ has a subalgebra isomorphic to
the three-dimensional Heisenberg algebra. Such a subalgebra cannot
be proper as it is not elementary. 
\par 
($\Leftarrow$) It is clear that the three-dimensional Heisenberg algebra 
is minimal non-elementary, so assume that $L$ is as described in (i). It
suffices to show that the maximal subalgebras of $L$ are
elementary. Let $M$ be a maximal subalgebra of $L$. Since Asoc $L
\leq M$ either $M = L^2$ or $M =$ Asoc $L + Fx$ for some $x \in L
\setminus L^2$. In the former case $M$ is abelian and so
elementary. So assume that $M =$ Asoc $L + Fx$ and $L = L^2 + Fx$.
Let $A$ be a minimal ideal of $M$. Then $A \leq$ Asoc $L$ and
$[A,L] \leq [A,x] \leq A$ so $A$ is an ideal of $L$. It follows
that Asoc $M =$ Asoc $L$ and $M$ splits over Asoc $M$. Hence $M$
is $\phi$-free. Since $M$ is clearly an $E$-algebra, it is
elementary. \bigskip

{\bf Note}: The algebras $L$ described in Theorem \ref{t:alel} (i)
are $A$-algebras (as every nilpotent subalgebra of $L$ is inside
$L^2$) that are not elementary. They are therefore not almost
algebraic.

\end{document}